\documentclass[11pt,a4paper]{amsart}

\usepackage[utf8]{inputenc}
\usepackage[T1]{fontenc}
\usepackage[english]{babel}

\usepackage{amsmath}
\usepackage{amssymb}
\usepackage{amsfonts}
\usepackage{mathtools}
\usepackage{mathrsfs}
\usepackage{amsthm}

\usepackage{graphicx}
\usepackage{tikz-cd}

\usepackage[colorlinks=true,
            linkcolor=blue,
            citecolor=blue,
            urlcolor=blue]{hyperref}

\theoremstyle{plain}
\newtheorem{theorem}{Theorem}[section]

\newtheorem{lemma}[theorem]{Lemma}
\newtheorem{corollary}[theorem]{Corollary}

\theoremstyle{definition}
\newtheorem{definition}[theorem]{Definition}

\theoremstyle{remark}
\newtheorem{remark}[theorem]{Remark}

\title[Semialgebraic Finiteness]
{Finiteness theorems for polynomials with complex coefficients}

\author{Samuel Martins}
\address{Departamento de Matemática,
Universidade Federal do Ceará,
Fortaleza, CE, Brazil}
\email{samuelsh9789@gmail.com}

\thanks{This study was financed in part by the
Coordenação de Aperfeiçoamento de Pessoal de
Nível Superior -- Brasil (CAPES) --
Finance Code 001.}

\begin{document}

\begin{abstract}
For fixed positive integers $n$ and $d$,
we prove that the space of complex polynomial
functions in $n$ variables, of degree at most
$d$ and vanishing at the origin, has only
finitely many local semialgebraic
$\mathcal R$-equivalence classes.
More precisely, the coefficient space admits
a finite partition into smooth, connected,
Zariski locally closed algebraic subsets
on which the equivalence type is constant.
\end{abstract}



\maketitle

% ==================================================
\section{Introduction}
% ==================================================

Let $n,m,d$ be positive integers and
$\mathbb K\in\{\mathbb R,\mathbb C\}$.
We denote by $P_0(n,m,d;\mathbb K)$
the space of polynomial mappings
$f:\mathbb K^n\to\mathbb K^m$
of degree at most $d$ satisfying $f(0)=0$.
This space can be identified with
$\mathbb K^l$ for some positive integer $l$.
For each $t\in\mathbb K^l$, we denote
the corresponding polynomial mapping by $f_t$.

Two polynomial mappings $f,g\in P_0(n,m,d;\mathbb{K})$
are said to be \emph{locally topologically
$\mathcal{R}$-equivalent} if there exist neighborhoods
$U,V$ of the origin in $\mathbb{K}^n$ and a
homeomorphism $h:U\to V$, with $h(0)=0$, such that
\[
f=g\circ h.
\]

Similarly, $f$ and $g$ are said to be
\emph{locally topologically
$\mathcal{RL}$-equivalent} if there exist
neighborhoods $U,V$ of the origin in
$\mathbb{K}^n$, neighborhoods $W_1,W_2$ of
the origin in $\mathbb{K}^m$, and homeomorphisms
$h:U\to V$ and $k:W_1\to W_2$, both fixing
the origin, such that
\[
k\circ f=g\circ h
\]
on a sufficiently small neighborhood of the origin.

Naturally, $\mathcal{R}$ and $\mathcal{RL}$ induce equivalence relations on $P_0(n,m,d; \mathbb{K})$. Fukuda \cite{Fukuda1976Types} proved that the number of equivalence classes under $\mathcal{RL}$-equivalence is finite when $m=1$. He later stated a local finiteness result for
topological $\mathcal{R}$-equivalence in the
real case \cite{Fukuda1978}, and a proof was
subsequently given by Koike and P\u{a}unescu
\cite{KoikePaunescu2022}.
One can, of course, define other analogous equivalence relations by imposing specific regularity conditions on the homeomorphisms $h$ (or the pair $h, k$), for instance, requiring them to be diffeomorphisms, semialgebraic homeomorphisms, bi-Hölder  or bilipschitz homeomorphisms.
In this article, we will be interested in the semialgebraic case.

Two polynomial mappings $f,g\in P_0(n,m,d;\mathbb{K})$
are said to be \emph{locally semialgebraically
$\mathcal{R}$-equivalent} if there exist neighborhoods
$U,V$ of the origin in $\mathbb{K}^n$ and a
semialgebraic homeomorphism $h:U\to V$,
with $h(0)=0$, such that
\[
f=g\circ h.
\]
Here, semialgebraicity is understood after identifying
$\mathbb{C}^n$ with $\mathbb{R}^{2n}$ when
$\mathbb{K}=\mathbb{C}$. Analogously, we can define local semi-algebraic $\mathcal{RL}$-equivalence.

Benedetti and Shiota \cite{BenedettiShiota1991}
proved the finiteness of semialgebraic
$\mathcal{RL}$-equivalence classes of real polynomial
functions of bounded degree. More recently, Koike
and P\u{a}unescu \cite{KoikePaunescu2022} proved that
the number of local semialgebraic
$\mathcal{R}$-equivalence classes in
$P_0(n,1,d;\mathbb{R})$ is finite.

The goal of this paper is to show that the number
of local semialgebraic $\mathcal{R}$-equivalence
classes in $P_0(n,1,d;\mathbb{C})$ is finite.
For this, we use the arc-wise analytic
trivialization theorem of Parusi\'nski and
P\u{a}unescu \cite{ParusinskiPaunescu2017}.

Moreover, we show that the coefficient space
admits a finite partition into smooth, connected,
Zariski locally closed complex algebraic
subvarieties on which the equivalence type is
constant. More precisely, we prove the following
theorem.

\begin{theorem}
For any positive integers $n$ and $d$, the number of
local semialgebraic $\mathcal{R}$-equivalence classes
in $P_0(n,1,d;\mathbb{C})$ is finite.

More precisely, there exists a finite partition
\[
\mathbb{C}^l=S_1\sqcup\cdots\sqcup S_r
\]
into smooth, connected, Zariski locally closed
complex algebraic subvarieties such that $f_s$
and $f_t$ are locally semialgebraically
$\mathcal{R}$-equivalent whenever $s,t\in S_j$.

\end{theorem}

Note that, although $\mathbb C$ has real dimension
two, the analogous finiteness result does not hold
in $P_0(n,2,d;\mathbb R)$ for any positive integers $n$ and $d$.
Indeed, the linear mappings $F_a(x_1,\ldots,x_n)=(x_1,ax_1),
 a\in\mathbb R$, are pairwise locally topologically
$\mathcal R$-inequivalent, since their image
germs at the origin are distinct.

The proof of Theorem~1.1 relies on the
arc-wise analytic trivialization theorem
of Parusi\'nski and P\u{a}unescu, stated
in Section~2. In Section~3, we construct
a suitable system of Weierstrass polynomials
satisfying the hypotheses of this theorem.
In Section~4, we establish the semialgebraic
properties needed to ensure that the
resulting trivializations are semialgebraic
for each fixed parameter.
In Section~5, we use these trivializations
to obtain local semialgebraic
$\mathcal R$-equivalences and prove
the main theorem by Noetherian induction.
Finally, Section~6 presents some
consequences of our results.

% ==================================================
\section{Preliminaries}
% ==================================================

Throughout the paper, we set $N=n+1$ and write
\[
x=(x_1,\ldots,x_N)=(y,z_1,\ldots,z_n),
\qquad
x^{(i)}=(x_1,\ldots,x_i).
\]

Let $W_i(t,x^{(i)})$ be a Weierstrass polynomial
in $x_i$, with coefficients holomorphic in
$(t,x^{(i-1)})$. We allow $W_i\equiv 1$ and
denote by $W_{i,\mathrm{red}}$ its reduced
Weierstrass polynomial. We use the convention
$\operatorname{Disc}(1)=1$.

We shall use the following consequence of
\cite[Theorem 3.3 and Proposition 3.7]
{ParusinskiPaunescu2017}.

\begin{theorem}
\label{thm:parusinski-paunescu}
Let $W_i(t,x^{(i)}), 1 \leq i \leq N$ be Weierstrass polynomials in $x_i$, possibly identically equal to $1$, with coefficients holomorphic in $(t,x^{(i-1)}), t \in (\mathbb C^r ,0)$.
Assume that $W_i\equiv 1$ implies
$W_j\equiv 1$ for every $j<i$,
$\operatorname{Disc}_{x_i}(W_{i,\mathrm{red}})
\mid x_1^{q_i}W_{i-1}, 2\leq i\leq N$, for some integers $q_i\geq 0$ and that $W_i(t,0, \cdots ,0,x_i)=x_{i}^{d_i}$ where $d_i=deg_{x_i}W_i$. Then
there exists an arc-wise analytic $\Phi(t,x)=(t,\phi_t(x))$ defined in a neighborhood of the origin, with the following properties:
\begin{enumerate}
\item $\phi_0=\mathrm{id}$ and
$\phi_t(0)=0$;

\item $\phi_t$ is a subanalytic
homeomorphism whose inverse is
arc-analytic;

\item $\phi_t$ has triangular form
\[
\phi_t(x)=
\bigl(x_1,\Psi_2(t,x^{(2)}),\ldots,
\Psi_N(t,x^{(N)})\bigr);
\]

\item For every $1\leq i\leq N$,
the induced map
\[
\phi_{i,t}(x^{(i)})
=
\bigl(x_1,\Psi_2(t,x^{(2)}),\ldots,
\Psi_i(t,x^{(i)})\bigr)
\]
satisfies
\[
\phi_{i,t}
\bigl(\{W_i(0,x^{(i)})=0\}\bigr)
=
\{W_i(t,x^{(i)})=0\}
\]
as germs at the origin.
\end{enumerate}
Moreover, the trivialization can be
constructed inductively by Whitney
interpolation of the roots of the
polynomials $W_i$.
\end{theorem}
\begin{proof}
For each $1\leq i\leq N$, set
\[
\widehat W_i=
\begin{cases}
W_i, & \text{if } W_i\not\equiv 1,\\
x_i, & \text{if } W_i\equiv 1.
\end{cases}
\]
Then each $\widehat W_i$ is a monic Weierstrass
polynomial of positive degree in $x_i$.
Moreover, by hypothesis, $\operatorname{Disc}_{x_i}
(\widehat W_{i,\mathrm{red}})
\mid x_1^{q_i}\widehat W_{i-1}, 2\leq i\leq N$. Indeed, if $W_i\equiv 1$, the discriminant
of $\widehat W_i=x_i$ is equal to $1$.

Set $F_0=1$, $F_N=\widehat W_N$, and $F_i=x_1^{q_{i+1}}\widehat W_i, 1\leq i<N$. Since
$\operatorname{Disc}_{x_1}
(\widehat W_{1,\mathrm{red}})=1$,
the system $\{F_i\}_{i=0}^N$ satisfies the
hypotheses of
\cite[Proposition 3.7]{ParusinskiPaunescu2017}.

Properties (1)--(3) follow from
\cite[Theorem 3.3 and Proposition 3.7]
{ParusinskiPaunescu2017}.
Property (4) follows from the regularity of
each $\widehat W_i$ with respect to the induced
trivialization whenever $W_i\not\equiv 1$,
and is immediate when $W_i\equiv 1$.

The construction by Whitney interpolation
is described in the proof of Theorem 3.3
of the same paper. In particular, whenever
$W_i\equiv 1$, the interpolation associated
with $\widehat W_i=x_i$ is the identity
in the $x_i$-coordinate.
\end{proof}
For further details on Whitney interpolation, we refer
the reader to Appendix A of
\cite{ParusinskiPaunescu2017}, and to Section 3
of the same paper for its use in the construction
of the functions $\Psi_i$.

\section{Algebraic Preparation and Weierstrass Systems}

Let $Y\subset\mathbb C^l$ be a closed irreducible
algebraic subvariety. For a Zariski open subset $Y^\circ\subset Y$,
we denote by $\mathbb{C}[Y^\circ]$ the ring
of regular functions on $Y^\circ$. More precisely,
a function $f:Y^\circ\to\mathbb{C}$ is regular
if, for every $p\in Y^\circ$, there exist a
Zariski open neighborhood $U\subset Y^\circ$
of $p$ and polynomials
$P,Q\in\mathbb{C}[t_1,\ldots,t_l]$,
with $Q$ nowhere zero on $U$, such that
\[
f|_U=\frac{P}{Q}\bigg|_U.
\]

In particular, if $Y^\circ$ is a principal
Zariski open subset of $Y$, that is,
\[
Y^\circ=D(h)
=\{t\in Y:h(t)\neq 0\}
\]
for some nonzero $h\in\mathbb{C}[Y]$, then
\[
\mathbb{C}[Y^\circ]
=
\mathbb{C}[Y]_h
=
\left\{
\frac{a}{h^k}:
a\in\mathbb{C}[Y],\ k\geq 0
\right\}
\subset\mathbb{C}(Y).
\]
For more details, see chapter 3 of \cite{Milne2024}. 

\begin{lemma}\label{lem:global-discriminant-chain}
Let $Y\subset\mathbb{C}^l$ be a closed irreducible
algebraic subvariety, with $Y\neq\{0\}$. Then there exist a linear automorphism
$A\in GL_n(\mathbb{C})$, independent of $t$,
a nonempty principal Zariski open subset $Y^\circ=D(h)\subset Y^{\mathrm{reg}}, h\in\mathbb{C}[Y]\setminus\{0\}$, an integer $1\leq i_0<N$, and integers
$m_i\geq 1$ and $q_i\geq 0$,
for $i_0<i\leq N$, together with polynomials
$P_i\in\mathbb{C}[Y^\circ][x^{(i)}],
 i_0\leq i\leq N$, and, for $i_0<i\leq N, P_i^*\in\mathbb{C}[Y^\circ][x^{(i)}], 
C_i\in\mathbb{C}[Y^\circ][x^{(i-1)}]$, such that the following conditions hold.

\begin{enumerate}
\item[\textup{(C1)}]

$P_N(t,x)
=x_1-f_t\bigl(A(x_2,\ldots,x_N)\bigr)$.

\item[\textup{(C2)}]
For every $i_0<i\leq N$, $P_i(t,0,\ldots,0,x_i)
=x_i^{m_i}\beta_i(t,x_i)$, where $\beta_i\in\mathbb{C}[Y^\circ][x_i]$
satisfies $\beta_i(t,0)\neq 0
\quad\text{for every }t\in Y^\circ$.

\item[\textup{(C3)}]
For every $i_0<i\leq N$, there exist
integers $r_i\geq 1$, $e_{i,k}\geq 1$,
and polynomials $Q_{i,k}\in
\mathbb{C}[Y^\circ][x^{(i)}],
 1\leq k\leq r_i$, of positive degree in $x_i$, which are
irreducible and pairwise nonassociated in
$\mathbb{C}(Y)[x^{(i)}]$, such that $P_i=C_i\prod_{k=1}^{r_i}Q_{i,k}^{e_{i,k}},
P_i^*=\prod_{k=1}^{r_i}Q_{i,k}$,
with $C_i(t,0)\neq 0
\quad\text{for every }t\in Y^\circ$.

\item[\textup{(C4)}]
For every $i_0<i\leq N$,
$\operatorname{Disc}_{x_i}(P_i^*)
=x_1^{q_i}P_{i-1}\not\equiv 0,
x_1\nmid P_{i-1}$.

\item[\textup{(C5)}]

$P_{i_0}(t,0)\neq 0
\quad\text{for every }t\in Y^\circ$.

\end{enumerate}
\end{lemma}

\begin{proof}
   
Let $R=\mathbb{C}[Y]$ be the ring of polynomials
in $l$ variables restricted to $Y$, and let
$K=\operatorname{Frac}(R)$.

We first construct, over $K$, the polynomials
$P_i$, $P_i^*$, and $C_i$, together with a linear
automorphism $A\in GL_n(\mathbb C)$, satisfying
(C1) and the following conditions:

\medskip
\noindent
$(C2')$ For every $i_0<i\leq N$,
$P_i(0,\ldots,0,x_i)=x_i^{m_i}\beta_i(x_i)$,
where $m_i\geq 1$, $\beta_i\in K[x_i]$,
and $\beta_i(0)\neq 0$ in $K$.

\medskip
\noindent
$(C3')$ For every $i_0<i\leq N$, there exist
$C_i\in K[x_1,\ldots,x_{i-1}]$ and
irreducible, pairwise nonassociated polynomials
$Q_{i,j}\in K[x_1,\ldots,x_i]$ of positive
degree in $x_i$ such that
$P_i=C_i\prod_{j=1}^{r_i}Q_{i,j}^{e_{i,j}}$
and $P_i^*=\prod_{j=1}^{r_i}Q_{i,j}$,
where $e_{i,j}\geq 1$ and
$C_i(0)\neq 0$ in $K$.

\medskip
\noindent
$(C4')$ For every $i_0<i\leq N$,
$\operatorname{Disc}_{x_i}(P_i^*)
=x_1^{q_i}P_{i-1}\neq 0$,
where $q_i\geq 0$ and $x_1\nmid P_{i-1}$.

\medskip
\noindent
$(C5')$ We have $P_{i_0}(0)\neq 0$ in $K$.

\medskip
We construct these polynomials by descending
induction on $i$. At step $i$, we assume that
$(I_i)$ holds, that is,
$0\neq P_i\in K[x_1,\ldots,x_i]$ and
$x_1\nmid P_i$, and that the polynomials
$P_k$, for $k>i$, already satisfy
(C2')--(C4').

The construction terminates at the first index
$i_0$ for which $P_{i_0}(0)\neq 0$, yielding
(C5').

    First, define $P_N = x_1 - f_t(x_2, \dots, x_N)$, where $f_t$ is non-zero in $K[x_2, \cdots , x_N]$ since $Y \neq \{0\}$. Thus, $P_N(0, x_2, \dots, x_N) = -f_t(x_2, \dots, x_N) \not\equiv 0$; that is, $x_1 \nmid P_N$. Let $H(x_2,\ldots,x_N)=P_N(0,x_2,\ldots,x_N)$,
and let $H_k$ be its first nonzero homogeneous
component. Since $H_k$ is homogeneous and nonzero and
$\mathbb C\subset K$ is an infinite subfield,
we can choose $\lambda_2,\ldots,\lambda_{N-1}\in\mathbb C$
such that
\[
H_k(\lambda_2,\ldots,\lambda_{N-1},1)\neq 0.
\]
Consider the triangular linear transformation
\[
T_N(x_1,\ldots,x_N)
=
(x_1,x_2+\lambda_2x_N,\ldots,
x_{N-1}+\lambda_{N-1}x_N,x_N).
\]
Then
\[
\begin{aligned}
&(P_N\circ T_N)(0,\ldots,0,x_N)\\
&=x_N^kH_k(\lambda_2,\ldots,\lambda_{N-1},1)
+O(x_N^{k+1}).
\end{aligned}
\]
Hence, after replacing $P_N$ by $P_N\circ T_N$,
we obtain
\[
P_N(0,\ldots,0,x_N)
=x_N^{m_N}\beta_N(x_N),
\quad m_N=k\geq 1,\quad \beta_N(0)\neq 0,
\]
where $\beta_N\in K[x_N]$.  This proves $C2'$.

Next, factoring $P_N$ in  $K[x_1,\ldots,x_{N-1}][x_N]$,
we obtain
\[
P_N=C_N\prod_{j=1}^{r_N}Q_{N,j}^{e_{N,j}},
\qquad
P_N^*=\prod_{j=1}^{r_N}Q_{N,j},
\]
where $C_N \in K[x_1,\ldots,x_{N-1}]$, is the greatest common divisor of the coefficients of  $P_N$  and the $Q_{N,j} \in K[x_1, \cdots , x_{N-1}][x_N]$
are pairwise non-associated irreducible factors of positive
degree in $x_N$. By ($C2'$), $C_N(0)\neq 0$, proving ($C3'$).

By Gauss's Lemma, $P_N^*$ is the product of distinct, non-associated irreducible factors in $K(x_1, \dots, x_{N-1})[x_N]$, where $K(x_1, \dots, x_{N-1})$ is the field of fractions of $K[x_1, \dots, x_{N-1}]$. Given that every irreducible polynomial over a field of characteristic zero is separable, it follows that $P_N^*$ is separable in $K(x_1, \dots, x_{N-1})[x_N]$. Consequently, the discriminant $Disc_{x_N}(P_N^*)$ is non-zero. We may write
\[
\operatorname{Disc}_{x_N}(P_N^*)
=x_1^{q_N}P_{N-1},
\qquad q_N\geq 0,\quad x_1\nmid P_{N-1},
\]
which proves ($C4'$) and initiates the descending induction.

Suppose $P_i$ satisfies $I_i$ and let $a_i = P_i(0) \in K$. If $P_i(0) \neq 0$, the process terminates. Otherwise, $a_i = 0$; in this case, it must be that $i > 1$, because if $i = 1$ and $P_1(0) = 0$, then $x_1 \mid P_1$, which contradicts $I_1$.
As in the initial step for $P_N$, there exists
a triangular linear transformation $T_i$,
fixing $x_1$ and all $x_j$ for $j\geq i$,
such that $P_i\circ T_i$ satisfies ($C2'$).

It remains to verify that $T_i$ preserves the
properties established at the previous steps
of the descending induction. For every $j>i$, the transformation $T_i$ fixes
$x_j$ and only changes variables preceding $x_j$.
In particular,
\[
(P_j\circ T_i)(0,\ldots,0,x_j)
=
P_j(0,\ldots,0,x_j),
\]
so the $x_j$-regularity of $P_j$ is preserved,
and hence ($C2'$) remains valid.

Moreover, $T_i$ induces an automorphism of
$K[x_1,\ldots,x_j]$ that preserves the subring
$K[x_1,\ldots,x_{j-1}]$. Consequently, the
factorization in ($C3'$) is preserved.

Since $T_i$ fixes $x_j$ and $x_1$, we also have
\[
\begin{aligned}
\operatorname{Disc}_{x_j}(P_j^*\circ T_i)
&=
\bigl(\operatorname{Disc}_{x_j}(P_j^*)\bigr)
\circ T_i\\
&=
x_1^{q_j}(P_{j-1}\circ T_i).
\end{aligned}
\]
Thus, ($C4'$) is preserved as well.

Finally, since $T_i$ fixes $x_1$, the condition
$x_1\nmid P_{j-1}$ remains valid.

As in the initial step for $P_N$, we define
$C_i$ as the greatest common divisor of the
coefficients of $P_i$ with respect to $x_i$.
Since (C2') holds, we have $C_i(0)\neq 0$
in $K$. Factoring the primitive part
$P_i/C_i$ proves (C3').

Again, by Gauss's lemma and characteristic zero,
$P_i^*$ is separable over
$K(x_1,\ldots,x_{i-1})$. Hence
$D_i=\operatorname{Disc}_{x_i}(P_i^*)\neq 0$.
Writing $D_i=x_1^{q_i}P_{i-1}$, with
$x_1\nmid P_{i-1}$, proves (C4')
and establishes $(I_{i-1})$.

The process repeats until we reach an $i_0$ such that $P_{i_0}(0) \neq 0$; this process is guaranteed to terminate, since for $i=1$, given that $x_1 \nmid P_1$, we have $P_1(0) \neq 0$.

As for (C1), we define $A\in GL_n(\mathbb C)$ by
\[
(T_N\circ\cdots\circ T_{i_0+1})(x_1,z)
=(x_1,Az),
\]
where $z=(x_2,\ldots,x_N)$.
This is possible since each $T_i$ fixes $x_1$
and acts only on the variables $x_2,\ldots,x_N$.
Thus, we obtain $P_N(t,x)=x_1-f_t(Az)$.

We now transfer these properties to $\mathbb{C}[Y^\circ]$ with original (C1)-(C5). Before let us clarify the differences between (C2')
and (C2). The condition (C2') states that
$\beta_i(0)=b_0\in K$ is nonzero, meaning
that $b_0$ is not identically zero as a
rational function on $Y$. At any point
$t\in Y$ where $b_0$ is defined,
$b_0(t)$ is an element of $\mathbb C$.
However, it is possible that $b_0\neq 0$
in $K$ while $b_0(t)=0$ for some such $t$.
We wish to find a nonempty Zariski open
subset $Y^\circ\subset Y$ such that $b_0$
is defined and $b_0(t)\neq 0$ for every
$t\in Y^\circ$. The same distinction applies
to the nonvanishing conditions in (C3')
and (C5').

There are only finitely many coefficients involved
in $P_i$, $P_i^*$, $C_i$, $Q_{i,k}$, and $\beta_i$.
Each coefficient can be written as $a/b$, where
$a,b\in R$, $b\neq 0$, and
$a(t)=\sum_{|\alpha|\leq d_1}a_\alpha t^\alpha$,
$b(t)=\sum_{|\beta|\leq d_2}b_\beta t^\beta$,
with $a_\alpha,b_\beta\in\mathbb C$ and
$t=(t_1,\ldots,t_l)\in Y$. Let $h_1 \in R \setminus \{0\}$ be the product of all these denominators, and consider
$$R_{h_1}=\{ \frac{a}{h_{1}^{m}}; a \in R, m \geq 0\}.$$

Note that $R \subset R_{h_1} \subset K$; consequently, identities that hold over $K$ remain true over $R_{h_1}$. Although these rational functions are nonzero
as elements of $K$, they may vanish at some
points of $D(h_1)$. Indeed, if $g=a/b\neq 0$
in $K$, with $a,b\in R$, then $a\neq 0$ in $R$,
meaning that $a$ is not identically zero on $Y$.
However, there may exist $t\in Y$ such that
$a(t)=0$, and hence $g(t)=0$ whenever
$b(t)\neq 0$.  Let $g_1, \dots, g_m \in R_{h_1} \setminus \{0\}$ be the collection consisting of $\beta_i(0)$, $C_i(0)$, and $P_{i_0}(0)$, such that we may write $g_j=a_j/h_{1}^{n_j}$, where $a_j \neq 0$ in $R$. Define $h_2=a_1 \cdots a_m$.  Finally, let $h_3 \in R$ be a non-zero function such that $h_3(t)=0$ for all $t \in Y^{sing}$; the desired set is then given by $Y^\circ=\{t \in Y; h_1(t)h_2(t)h_3(t) \neq 0 \}$.

Since $Y$ is irreducible, $R$ is an integral domain, and hence
$h_1h_2h_3 \neq 0$. Therefore,
\[
Y^\circ = \{t \in Y : h_1(t)h_2(t)h_3(t) \neq 0\}
\]
is a nonempty Zariski-open subset of $Y^{\mathrm{reg}}$.
By construction, all coefficients belong to
$\mathbb C[Y^\circ]=R_h$, where $h=h_1h_2h_3$.
Since $R_h\subset K$, all polynomial identities
established over $K$ remain valid over $R_h$.

Moreover, the choice of $h_2$ guarantees that
\[
\beta_i(t,0)\neq 0,\qquad
C_i(t,0)\neq 0,\qquad
P_{i_0}(t,0)\neq 0
\]
for every $t\in Y^\circ$.
Finally, the condition $x_1\nmid P_{i-1}$
is preserved, since
\[
R_h[x_1,\ldots,x_{i-1}]
\hookrightarrow
K[x_1,\ldots,x_{i-1}]
\]
is injective.

Thus, all properties (C1)--(C5) hold on
$Y^\circ$, completing the proof.

\end{proof}
We will now construct the system of Weierstrass polynomials that satisfies the hypotheses of Theorem 2.1.
\begin{lemma}
Let $Y^\circ$, $A$, $i_0$, and $P_i$ be as in
Lemma~3.1. Fix $t_0\in Y^\circ$ and set
$r=\dim_{\mathbb C}Y$. Choose local holomorphic
coordinates $u\in\mathbb C^r$ on $Y^\circ$
centered at $t_0$, and denote by $t=t(u)$
the corresponding local parametrization,
with $t(0)=t_0$. Set $W_i\equiv 1$ for $i\leq i_0$. For $i>i_0$,
let $W_i$ be the Weierstrass polynomial associated
with $P_i$ with respect to $x_i$ at $(t_0,0)$.
Then the following properties hold:

\begin{enumerate}
\item[\textup{(a)}]
For every $i>i_0$, we have $W_i(u,0,\ldots,0,x_i)=x_i^{m_i}$.

\item[\textup{(b)}]
For every $i>i_0$, $\operatorname{Disc}_{x_i}(W_{i,\mathrm{red}})
\mid x_1^{q_i}W_{i-1}
\quad\text{in }\mathcal O_{\mathbb C^{r+i-1},0}$, where $W_{i_0}=1$.

\item[\textup{(c)}]
In a neighborhood of the origin of $\mathbb{C}^r \times \mathbb{C}^N$ $
\{W_N(u,x)=0\}
=
\{x_1=f_{t(u)}(A(x_2,\ldots,x_N))\}$.

\end{enumerate}
\end{lemma}

\begin{proof}
   For simplicity, set $\mathcal A_{i-1}
=\mathcal O_{\mathbb C^{r+i-1},0}$,
the ring of holomorphic germs at the origin
in the variables $(u,x^{(i-1)})$.For $i > i_0$, by (C2), $P_i$ is $x_i$-regular of order $m_i$ at $(0,0)$. By the Weierstrass Preparation Theorem, there exists a unit $U_i$ and a Weierstrass polynomial $W_i$ such that $P_i = U_i W_i$. Evaluating this at $x^{(i-1)} = 0$ and comparing it with $P(u, 0, \dots, 0, x_i) = x_i^{m_i} \beta_i(u, x_i)$ yields $W_i(u, 0, x_i) = x_i^{m_i}$. This proves (a).

    By (C3), $P_i$ and $P_i^*$ define the same germ of a zero set. Applying the Weierstrass Preparation Theorem, we can write $P_i^* = V_i H_i$, where $V_i$ is a unit and $H_i$ is a Weierstrass polynomial. Note that $H_i$ and $W_i$ define the same zero set. We wish to show that $H_i$ is reduced.

Since $P_i^*$ is a polynomial in $x_i$ with coefficients in $\mathcal{A}_{i-1}$ and $H_i$ is monic in $x_i$, we can perform division within $\mathcal{A}_{i-1}[x_i]$; by the uniqueness of division, we obtain $P_i^* = H_i R_i$, where $R_i \in \mathcal{A}_{i-1}[x_i]$ and $R_i = V_i$ is a unit.

By the discriminant product formula and the
homogeneity of the discriminant
(see \cite[Theorems 3.10 and 3.5(i)]{JansonDiscriminant}),
we have that
$\operatorname{Disc}_{x_i}(H_i)$ divides
$\operatorname{Disc}_{x_i}(P_i^*)$
in $\mathcal{A}_{i-1}$.
By (C4), $\operatorname{Disc}_{x_i}(P_i^*)\neq 0$. Indeed, if $\operatorname{Disc}_{x_i}(P_i^*)$
were identically zero as a holomorphic germ
in $(u,x^{(i-1)})$, then all its coefficients,
regarded as regular functions on $Y^\circ$,
would vanish on a Euclidean neighborhood
of $t_0$. Since $Y^\circ$ is a connected
complex manifold, the identity theorem
would imply that all these coefficients
vanish identically on $Y^\circ$,
contradicting (C4).
Hence $\operatorname{Disc}_{x_i}(H_i)\neq 0$,
so $H_i$ is reduced. Since $P_i$ and $P_i^*$
define the same germ of zero set, the same
holds for $W_i$ and $H_i$. By uniqueness
of the reduced Weierstrass polynomial,
we obtain $H_i=W_{i,\mathrm{red}}$.

Thus, $\operatorname{Disc}_{x_i}(W_{i,\mathrm{red}})
\mid
\operatorname{Disc}_{x_i}(P_i^*)
=
x_1^{q_i}P_{i-1}$.

If $i-1>i_0$, then $P_{i-1}=U_{i-1}W_{i-1}$,
where $U_{i-1}$ is a unit. Hence, $\operatorname{Disc}_{x_i}(W_{i,\mathrm{red}})
\mid x_1^{q_i}W_{i-1}$.

If $i-1=i_0$, property (C5) implies that
$P_{i_0}$ is a unit at the point under consideration.
Therefore, the divisibility reduces to
\[
\operatorname{Disc}_{x_i}(W_{i,\mathrm{red}})
\mid x_1^{q_i}
=
x_1^{q_i}W_{i_0},
\]
since $W_{i_0}=1$.

Finally, by (C1),
\[
P_N(u,x)=x_1-f_{t(u)}(A(x_2,\ldots,x_N)).
\]
Since $P_N=U_NW_N$, with $U_N$ a holomorphic
unit, $P_N$ and $W_N$ define the same germ
of zero set. This proves (c).

\end{proof}
\section{Semialgebraicity of the Trivialization  }

As a consequence of the next theorem, we will have that for each fixed $u$ sufficiently
close to the origin, the coefficients of
$W_i(u,x^{(i)})$ from Lemma 3.2 are Nash
germs in $x^{(i-1)}$.

\begin{theorem}
Let $P(a,w)\in\mathbb{C}[a_1,\ldots,a_r,w]$
be such that $P(0,w)$ has order $m\geq 1$
at $w=0$. Let $W(a,w)=w^m+\sum_{j=1}^{m}c_j(a)w^{m-j}$ be the Weierstrass polynomial associated with
$P$ with respect to $w$ at the origin. Then there exists a semialgebraic neighborhood
$U\subset\mathbb{C}^r$ of the origin such that
each coefficient function $c_j:U\to\mathbb{C}$
is semialgebraic.

\end{theorem}

\begin{proof}
    We can write $P(0,w)=w^m \beta(w)$ with $\beta(0) \neq 0$. By the Weierstrass Preparation Theorem, we have $P(a,w)=U(a,w) W(a,w)$ in a neighborhood $0 \in \Omega \subset \mathbb{C}^r \times \mathbb{C}$, where $U$ is a unit and $W$ is a Weierstrass polynomial. 
Choose $\rho>0$ and $\epsilon_0>0$ such that
$P(0,w)\neq 0$ for $0<|w|\leq 2\rho$ and
$B_{\epsilon_0}\times\{|w|\leq 2\rho\}\subset\Omega$.
 It follows that   $P(0, \cdot)$ has no roots  in $|w|= \rho$ and it has exactly $m$ roots inside, counting multiplicity. 
Since $P$ is continuous, there exists
$0<\epsilon<\epsilon_0$ such that, for
$B=\{a\in\mathbb{C}^r:|a|^2<\epsilon^2\}$,
\[
|P(a,w)-P(0,w)|<|P(0,w)|
\]
for every $a\in B$ and $|w|=\rho$.

    Using Rouché's Theorem for $P(0,w)$ and $P(a,w)-P(0,w)$, it follows that for every $a \in B$, the number of zeros of $P(a, \cdot)$ in the disk $\{|w|< \rho\}$, counted with multiplicity, is $m$.

We shall now prove that if $a \in B$ and $W(a,w)=0$, then $|w| < \rho$. Indeed, if $|w|=\rho$, we have $|c_1(a)w^{m-1}+ \cdots + c_m(a)| \leq \sum |c_j(a)| \rho^{m-j}$. Since $c_j(0)=0$, we may assume—by shrinking $B$ if necessary—that $|c_j(a)|<\frac{\rho^j}{2m}$, such that $|\sum c_j(a) w^{m-j}| < \rho^m$. Consequently, we can apply Rouché's Theorem to $w^m$ and $c_1(a)w^{m-1}+ \cdots + c_m(a)$; it follows that $W(a, \cdot)$ has exactly $m$ zeros (counted with multiplicity) in the disk $\{|w|< \rho\}$. Since $W(a, \cdot)$ has at most $m$ roots, the claim is proved.

Let $D = \deg_w P$. For $c=(c_1, \dots, c_m) \in \mathbb{C}^m$, we define the polynomial $Q_c(w) = w^m + c_1 w^{m-1} + \dots + c_m$.

\textbf{Claim:} for all $a \in B$, the polynomial $Q_c$ coincides with $W(a , \cdot)$ if and only if there exists a polynomial $R(w)=b_0+b_1w+ \cdots b_{D-m}w^{D-m} \in \mathbb{C}[w]$ that satisfies:

\begin{enumerate}
    \item $P(a,w)=Q_c(w) R(w)$;
    \item $Q_c(x)=0 \implies |x|< \rho$;
    \item $|x|  < \rho \implies R(x) \neq 0.$
\end{enumerate}

Suppose first that $Q_c=W(a,\cdot)$.
Condition (2) has already been established. We will now show 1. By performing polynomial division in $\mathcal{O}_{\mathbb{C}^r,0}[w]$, we obtain $R \in \mathcal{O}_{\mathbb{C}^r,0}[w]$ and $S$ such that $\deg_w S < m$, with $P = WR + S$. However, by the Weierstrass Division Theorem, we already have $P = WU + 0$. Thus, by uniqueness, it follows that $S = 0$ and $R = U$.  Since $P(a,\cdot)$ and $W(a,\cdot)$ have exactly
$m$ zeros in $\{|w|<\rho\}$, counted with
multiplicity, $R$ has no zeros in this disk,
which proves (3).

Conversely, suppose that (1)--(3) hold.
Then the zeros of $P(a,\cdot)$ in
$\{|w|<\rho\}$ are precisely those of
$Q_c$, counted with multiplicity.
The same holds for $W(a,\cdot)$.
Since $Q_c$ and $W(a,\cdot)$ are monic
polynomials of degree $m$, we obtain
$Q_c=W(a,\cdot)$.

Let
\[
\mathcal P_{D-m}
=
\{R\in\mathbb C[w]:\deg R\leq D-m\}.
\]
Note that conditions 1, 2, and 3 are semialgebraic; therefore, the following set is also semialgebraic:
\[
S=
\left\{
(a,c)\in B\times\mathbb C^m:
\begin{array}{l}
\text{there exists }R\in\mathcal P_{D-m}\\
\text{satisfying conditions (1)--(3)}
\end{array}
\right\}.
\]
By the statement, for every $a \in B$, there exists $c \in \mathbb{C}^m$ such that $(a,c) \in S$, with $c(a) = (c_1(a), \dots, c_m(a))$; that is, the graph of $c|_B$ is $S$, which proves the result.
\end{proof}

\begin{lemma}
Assume the hypotheses of Theorem 2.1, and let
$\Phi$ be the trivialization obtained by Whitney
interpolation. Fix $u$ sufficiently close to $0$
and suppose that the coefficients of
$W_i(0,\cdot)$ and $W_i(u,\cdot)$ are Nash germs. Then, after restricting to a sufficiently small
semialgebraic polydisc $\Delta\subset\mathbb{C}^N$
centered at the origin, the map $\Phi_u:\Delta\longrightarrow\mathbb{C}^N$ is semialgebraic.  
\end{lemma}

\begin{proof}
    To prove that $\Phi_u$ is semialgebraic, it suffices to show that each $\Psi_j$ is also semialgebraic. Let $\Delta^i = \Delta_1 \times \cdots \times \Delta_i$ be a polydisc in $\mathbb{C}^i$ such that the previously obtained identities hold, and such that $P_{s,x'} = W_i(su, \phi_{i-1,su}(x'), w)$ (for $s \in [0, 1]$) has no roots on its boundary, where  \newline
    $\phi_{i,u}(x^{(i)}) = (\Psi_1(u, x^{(1)}), \Psi_2(u, x^{(2)}), \dots, \Psi_i(u, x^{(i)}))$ and $\Psi_1(u,x^{(1)})=x_1$. The proof will proceed by induction on $\phi_i$. If $i=1$, we have $\phi_1(x_1) = x_1$, which is obviously semialgebraic. We assume that $\phi_{i-1}$ is semialgebraic and will show that $\phi_i$ is semialgebraic. If $W_i \equiv 1$, then $\phi_i = \phi_{i-1} \times \text{Id}$, which is semialgebraic; thus, we may assume $d_i \geq 1$.

    Define $E=\{x' \in \Delta^{i-1} ; x_1 W_{i-1}(0,x') \neq 0\}$. Since $W_{i-1}(0, \cdot) $ is semialgebraic, we have that $E$ is semialgebraic. Since $W_{i-1}$ is either a constant equal to 1 or a monic polynomial in $x_{i-1}$, it follows that $x_1W_{i-1}(0,x')$ is non-zero; thus, $E$ is dense in $\Delta^{i-1}$. In fact, for every $x' \in E$ and $0 \leq s \leq 1$, we have $x_1W_{i-1}(su, \phi_{i-1,su}(x')) \neq 0$. This is because, by the construction of $\phi_{i-1}$, the condition $x_1W_{i-1}(0,x')=0$ holds if and only if $x_1W_{i-1}(u,\phi_{i-1,u}(x'))=0$. Since $Disc_{x_i}(W_{i,red}) \mid x_{1}^{q_i} W_{i-1}$, it follows that $Disc_{x_i}(W_{i,red})$ does not vanish at $(su,\phi_{i-1,su}(x'))$ for any $s \in [0,1]$ and $x' \in E$.

    \textbf{Claim:} $W_i(su, \phi_{i-1,su}(x'),w)$ has exactly $\rho=deg_{x_i} W_{i,red}$  distinct roots for every $s \in [0,1]$ and for all $x' \in E$. 

    Indeed, let $x' \in E$ and $s \in [0,1]$, and define $Q_s(w) = W_{i,red}(su, \phi_{i-1,su}(x'),w)$. We can write $Q_s(w) = w^{\rho} + c_1(s) w^{\rho -1} + \cdots + c_{\rho}(s)$. Based on what was established previously, the discriminant $Disc(Q_s)$ is non-zero for all $s \in [0,1]$. This is equivalent to saying that all roots of $Q_s$ are simple; therefore, $Q_s$ has $\rho$ roots for every $s \in [0,1]$. This proves the claim.

    By the continuity of the roots of a polynomial, for each $s_0 \in [0,1]$, there exist continuous functions $\lambda_1, \dots, \lambda_{\rho}$ defined on an open neighborhood of $s_0$ such that $Q_s(\lambda_j(s)) = 0$. Since the discriminant of $Q_s$ does not vanish
on $[0,1]$, its roots can be continuously
enumerated along the entire interval. Furthermore, since $Q_s$ has only simple roots,
we have $\lambda_j(s)\neq\lambda_k(s)$
for every $s\in[0,1]$ and $j\neq k$.

    \textbf{Claim:} For each  $j$, there exists an integer $m_j$ with $mult_{\lambda_j(s)}P_{s,x'}=m_j$.

    Indeed, by writing $W_i = H_1^{e_1} \cdots H_r^{e_r}$ as a product of irreducibles, we have $W_{i,red} = H_1 \cdots H_r$. Since all roots are simple in $W_{i, red}$, every root of $W_i$ is a root of exactly one $H_l$, with multiplicity $e_l$. Consequently, as this root varies along $\lambda_j$, it must always remain a root of the same $H_l$; otherwise, at some point, $\lambda_j$ would be a root of two distinct irreducible factors. This proves the assertion.

For $x'\in E$, let
$\mu_{s,x'}(a)$ denote the multiplicity
of $a$ as a root of $P_{s,x'}$.
Define
\[
M=
\left\{
(x',A,B)\in E\times\mathbb C^{2\rho}:
\begin{array}{l}
A=(A_1,\ldots,A_\rho)
\text{ enumerates the distinct}\\
\text{roots of }P_{0,x'},\\
B=(B_1,\ldots,B_\rho)
\text{ enumerates the distinct}\\
\text{roots of }P_{1,x'},\\
\mu_{0,x'}(A_j)=\mu_{1,x'}(B_j),
\quad 1\leq j\leq\rho
\end{array}
\right\}.
\]
which is semialgebraic.

We claim that the projection
\[
\pi:M\longrightarrow E,\qquad
\pi(x',A,B)=x',
\]
is a finite covering map.

Indeed, fix $x'_0\in E$. Since the
discriminants of $P_{0,x'}^{\mathrm{red}}$
and $P_{1,x'}^{\mathrm{red}}$ do not vanish
on $E$, there exists a sufficiently small
neighborhood $U\subset E$ of $x'_0$
and continuous functions
\[
a_1,\ldots,a_\rho,\,
b_1,\ldots,b_\rho:U\longrightarrow\mathbb C
\]
enumerating the distinct roots of
$P_{0,x'}$ and $P_{1,x'}$, respectively.

The multiplicities of these roots are
constant on $U$. Thus, every admissible
pair of permutations $(\sigma,\tau)$
determines a continuous section
\[
F_{\sigma,\tau}(x')
=
\bigl(x',
(a_{\sigma(j)}(x'))_{j=1}^{\rho},
(b_{\tau(j)}(x'))_{j=1}^{\rho}\bigr),
\]
where admissibility means that the
corresponding multiplicities agree.

The images of these sections are pairwise
disjoint and
\[
\pi^{-1}(U)
=
\bigsqcup_{(\sigma,\tau)}
F_{\sigma,\tau}(U).
\]
Each restriction of $\pi$ to
$F_{\sigma,\tau}(U)$ is a homeomorphism
onto $U$. Hence $\pi$ is a finite
covering map.

Define

\[
M_{\mathrm{cont}}
=
\left\{
\begin{aligned}
(x',A,B)\in M:\;&
B_j \text{ is the endpoint of the root}\\
&\text{obtained by continuously transporting } A_j\\
&\text{along the path } s\mapsto su
\text{ in the family}\\
& w\mapsto W_i(su,\phi_{i-1,su}(x'),w),\\
&\hfill 0\leq s\leq 1,\quad 1\leq j\leq\rho
\end{aligned}
\right\}.
\]

We claim that $M_{\mathrm{cont}}$ is both open
and closed in $M$. Fix $x'_0\in E$.
Since the discriminant of
\[
P_{s,x'}^{\mathrm{red}}(w)
= W_{i,\mathrm{red}}
(su,\phi_{i-1,su}(x'),w)
\]
does not vanish at $(s,x'_0)$ for any
$s\in[0,1]$, there exists a sufficiently
small neighborhood $U\subset E$ of $x'_0$
on which the roots admit continuous
enumerations
\[
\lambda_1,\ldots,\lambda_\rho:
[0,1]\times U\longrightarrow\mathbb C.
\]
For each $x'\in U$, the correspondence
determined by continuation is
\[
\lambda_j(0,x')\longmapsto
\lambda_j(1,x'),\qquad 1\leq j\leq\rho.
\]
Thus $M_{\mathrm{cont}}\cap\pi^{-1}(U)$
is a union of sheets of the finite covering
$\pi:M\to E$. Hence $M_{\mathrm{cont}}$
is open and closed in $M$, and therefore
a union of connected components of $M$.
Since $M$ is semialgebraic, so is
$M_{\mathrm{cont}}$.

For each $(x',A,B)\in M_{\mathrm{cont}}$,
we associate $\bar A,\bar B\in\mathbb C^{d_i}$,
obtained by repeating $A_j$ and $B_j$ according
to their multiplicities. Since there are only
finitely many possible multiplicity patterns,
the map
\[
(x',A,B)\longmapsto(\bar A,\bar B)
\]
is semialgebraic on $M_{\mathrm{cont}}$.

Recall that
\[
\Psi_i(u,x',w)
=\psi_{d_i}(w;a(0,x'),a(1,x')),
\]
where $\psi_{d_i}$ is the Whitney interpolation,
$a(0,x')$ is the list of roots of
$W_i(0,x',\cdot)$ counted with multiplicity,
and $a(1,x')$ is the corresponding list of roots
of $W_i(u,\phi_{i-1,u}(x'),\cdot)$ obtained
by continuous transport.

Choosing the enumerations
$a(0,x')=\bar A$ and $a(1,x')=\bar B$,
we obtain
\[
\Psi_i(u,x',w)
=\psi_{d_i}(w;\bar A,\bar B).
\]

It follows from the explicit formula for Whitney
interpolation in \cite[Example A.7, formula (A.13)]
{ParusinskiPaunescu2017} that $\psi_{d_i}$ is
semialgebraic. Moreover, by \cite[Remark A.2]
{ParusinskiPaunescu2017}, for every permutation
$\sigma\in S_{d_i}$, we have
\[
\psi_{d_i}(w;\sigma(\bar A),\sigma(\bar B))
=
\psi_{d_i}(w;\bar A,\bar B),
\]

where
\[
\sigma(a_1,\ldots,a_{d_i})
=(a_{\sigma(1)},\ldots,a_{\sigma(d_i)}).
\]

In particular, the value of
$\psi_{d_i}(w;\bar A,\bar B)$ is independent
of the chosen enumeration of the roots.

Therefore, the following set is semialgebraic:
\[
H=
\left\{
(x',w,A,B,v):
\begin{array}{l}
x'\in E,\quad w\in\Delta_i,\\
(x',A,B)\in M_{\mathrm{cont}},\\
v=\psi_{d_i}(w;\bar A,\bar B)
\end{array}
\right\}.
\]
By applying the projection
\[
p(x',w,A,B,v)=(x',w,v),
\]
it follows from the Tarski--Seidenberg theorem
that the following set is also semialgebraic:
\[
G=
\left\{
(x',w,v):
\begin{array}{l}
x'\in E,\quad w\in\Delta_i,\\
\exists\,(A,B)\in M_{\mathrm{cont}}(x')\\
\text{such that }
v=\psi_{d_i}(w;\bar A,\bar B)
\end{array}
\right\}.
\]
For each $x'\in E$, set $M_{\mathrm{cont}}(x')
=
\{(A,B):(x',A,B)\in M_{\mathrm{cont}}\}$.

Since $graph(\Psi_i(u, \cdot )|_{E \times \Delta_i})=G$, it follows that $\Psi_i(u, \cdot )|_{E \times \Delta_i}$ is semialgebraic; then, using the fact that $E \times \Delta_i$ is dense in $\Delta^i$ and extending by continuity, we conclude that $\Psi_i(u, \cdot)$ is semialgebraic on $\Delta^i$. This completes the induction process and the proof.

\end{proof}
%==================================================
\section{Proof of the Main Theorem }

\begin{lemma}
    Let $Y$ and $Y^\circ$ be as in Lemma ~3.1, and let $t_0 \in Y^\circ$. Then there exists an open set $B$ in $Y^\circ$ containing $t_0$ such that $f_t$ and $f_{t_0}$ are locally semialgebraically $R$-equivalent for all $t \in B$.
\end{lemma}

\begin{proof}
Fix $t_0\in Y^\circ$ and choose local holomorphic
coordinates $u$ centered at $t_0$, as in
Lemma~3.2. Let $\Phi(u,x)=(u,\phi_u(x))$
be the trivialization given by Theorem~2.1.
For each $t$ sufficiently close to $t_0$, set
$\Phi_t=\phi_{u(t)}$.

Write $g_t(z)=f_t(Az)$ and $V_t=\{(y,z)\in\mathbb C\times\mathbb C^n:
y=g_t(z)\}$. By Lemma~3.2 and Theorem~2.1, we have
$\Phi_t(V_{t_0})=V_t$ as germs at the origin.

Consider the mappings $\gamma_{t_0}(z)=(g_{t_0}(z),z),
 \pi(y,z)=z$, and define
\[
\sigma_t=\pi\circ\Phi_t\circ\gamma_{t_0}.
\]
By Theorem~4.1 and Lemma~4.2, $\sigma_t$ is
a local semialgebraic homeomorphism fixing
the origin.

Since $\Phi_t$ preserves the first coordinate,
we have $\Phi_t(g_{t_0}(z),z)
=(g_{t_0}(z),\sigma_t(z))\in V_t$. Hence $g_t\circ\sigma_t=g_{t_0}$.

Setting $h_t=A\circ\sigma_t\circ A^{-1}$, we obtain
\[
f_t\circ h_t=f_{t_0}.
\]
This proves the result.
\end{proof}
\begin{lemma}
\label{lem:generic-equivalence}
Let $Y\subset\mathbb{C}^l$ be a closed  irreducible
algebraic subvariety, with $Y\neq\{0\}$.
Then the Zariski open subset $Y^\circ$ given by
Lemma~\ref{lem:global-discriminant-chain}
is entirely contained in a single local
semialgebraic $\mathcal{R}$-equivalence class.
\end{lemma}
\begin{proof}
By Lemma 5.1, the equivalence classes are
open in $Y^\circ$. Since $Y^\circ$ is
irreducible, it is connected in the Euclidean
topology. Hence $Y^\circ$ \textit{is contained in a single equivalence class} 
\end{proof}

\begin{proof}[Proof of Theorem~1.1]
We proceed by Noetherian induction on the
Zariski closed subsets of $\mathbb C^l$. For a closed algebraic subset
$Y\subset\mathbb C^l$, let $\mathcal P(Y)$
denote the assertion that $Y$ admits a finite
partition into smooth, connected, Zariski
locally closed complex algebraic subvarieties,
on each of which the local semialgebraic
$\mathcal R$-equivalence type is constant. The assertion is immediate for
$Y=\varnothing$ and $Y=\{0\}$.

Let $Y\subset\mathbb C^l$ be a nonempty
Zariski closed subset, with $Y\neq\{0\}$,
and assume that $\mathcal P(Z)$ holds
for every proper Zariski closed subset
$Z\subsetneq Y$. Write $Y=Y_1\cup\cdots\cup Y_r$ for the decomposition of $Y$ into its
distinct irreducible components, choosing
$Y_1\neq\{0\}$.

By Lemma~\ref{lem:generic-equivalence},
there exists a nonempty Zariski open subset
$Y_1^\circ\subset Y_1^{\mathrm{reg}}$ such that all polynomials $f_t$, with
$t\in Y_1^\circ$, are locally
semialgebraically $\mathcal R$-equivalent. Set
\[
T=\bigcup_{j=2}^r Y_j,
\qquad
S=Y_1^\circ\setminus T,
\]
where $T=\varnothing$ if $r=1$. Since $Y_1$ is an irreducible component
of $Y$, it is not contained in $T$.
Thus $Y_1\setminus T$ is a nonempty
Zariski open subset of $Y_1$.
As $Y_1$ is irreducible, we obtain
$S=Y_1^\circ\cap(Y_1\setminus T)
\neq\varnothing$.

Moreover, $S$ is a Zariski open subset
of the irreducible variety $Y_1$.
Hence $S$ is irreducible and therefore
connected in the Euclidean topology.
Since $S\subset Y_1^{\mathrm{reg}}$,
it is smooth. Furthermore, $S$ is
Zariski locally closed in $\mathbb C^l$,
because $Y_1$ is closed and $S$ is
open in $Y_1$.

By construction, $S\subset Y_1^\circ$,
so the local semialgebraic
$\mathcal R$-equivalence type is
constant on $S$. Now observe that
$Y\setminus S
=
(Y_1\setminus Y_1^\circ)\cup T$. Consequently, $Y\setminus S$ is a
proper Zariski closed subset of $Y$.
By the induction hypothesis, it admits
a finite partition $Y\setminus S
=
S_2\sqcup\cdots\sqcup S_k$ into smooth, connected, Zariski locally
closed complex algebraic subvarieties
on which the equivalence type is constant.

Therefore,
\[
Y=S\sqcup S_2\sqcup\cdots\sqcup S_k
\]
is a partition with the required
properties. This proves $\mathcal P(Y)$.

By Noetherian induction, $\mathcal P(Y)$
holds for every Zariski closed subset
$Y\subset\mathbb C^l$. In particular,
taking $Y=\mathbb C^l$ proves the theorem.
\end{proof}
\begin{remark}
In Theorem~1.1, if $s,t\in S_j$, then
$f_s$ and $f_t$ are locally $\mathcal R$-equivalent
via a semialgebraic bi-arc-analytic homeomorphism;
see \cite[Theorem 3.3]{ParusinskiPaunescu2017}.
\end{remark}
\section{Applications }

A subset $A\subseteq\mathbb C^l$ is called
Zariski constructible if it can be written as
a finite union
\[
A=\bigcup_{i=1}^k(B_i\setminus C_i),
\]
where $B_i$ and $C_i$ are algebraic subsets
of $\mathbb C^l$ (see \cite{Milne2024}).
In particular, every Zariski constructible set
is semialgebraic. Since each local semialgebraic
$\mathcal R$-equivalence class is a finite union
of some of the Zariski locally closed subsets
$S_j$ given by Theorem~1.1, we obtain the
following corollary.

\begin{corollary}
For any positive integers $n$ and $d$, each local
semialgebraic $\mathcal R$-equivalence class in
$P_0(n,1,d;\mathbb C)$ is a Zariski constructible
subset of the coefficient space $\mathbb C^l$
and, consequently, is semialgebraic.
\end{corollary}

\begin{corollary}
    In $P_0(n,1,d;\mathbb{C})$, there is a finite number of  local semialgebraic $\mathcal{RL}$ classes and  there is a finite number of local subanalytic $\mathcal{R}$ and $\mathcal{RL}$ classes.
\end{corollary}

\begin{corollary}
Let $\{f_t\}_{t\in\mathbb C}$ be a family of
polynomials in $P_0(n,1,d;\mathbb C)$ of the form
\[
f_t(x)=
\sum_{1\leq|\beta|\leq d}
a_\beta(t)x^\beta,
\qquad a_\beta\in\mathbb C[t].
\]
Then there exists a finite subset
$E\subset\mathbb C$ such that $f_s$ and $f_t$
are locally semialgebraically
$\mathcal R$-equivalent for all
$s,t\in\mathbb C\setminus E$.
\end{corollary}

\begin{proof}
Let $\alpha:\mathbb C\longrightarrow\mathbb C^l$, 
$\alpha(t)=(a_\beta(t))_{1\leq|\beta|\leq d}$, be the polynomial map parametrizing the family,
and set $Y=\overline{\alpha(\mathbb C)}^{\,\mathrm{Zar}}$. Since $\mathbb C$ is irreducible and $\alpha$
is continuous in the Zariski topology,
$Y$ is an irreducible algebraic variety.

If $Y$ is a point, the assertion is immediate.
Otherwise, by
Lemma~\ref{lem:generic-equivalence},
there exists a nonempty Zariski open subset
$Y^\circ\subset Y$ on which the local
semialgebraic $\mathcal R$-equivalence
type is constant.

Since $\alpha(\mathbb C)$ is Zariski dense
in $Y$, the set $\alpha^{-1}(Y^\circ)$
is nonempty. Moreover, it is Zariski open
in $\mathbb C$, because $\alpha$ is polynomial. Consequently, $E=\mathbb C\setminus\alpha^{-1}(Y^\circ)$ is a proper Zariski closed subset of
$\mathbb C$, and hence is finite. For every $s,t\in\mathbb C\setminus E$,
we have $\alpha(s),\alpha(t)\in Y^\circ$.
Therefore, $f_s$ and $f_t$ are locally
semialgebraically $\mathcal R$-equivalent.
\end{proof}

By the arc-wise analytic properties of the
trivialization constructed in
\cite[Theorem 3.3]{ParusinskiPaunescu2017},
the homeomorphisms $\Phi_u$ and
$\Phi_u^{-1}$ are arc-analytic.
 Consequently, the mapping $h$ such that $f=g \circ h$ is semialgebraic and bi-arc-analytic, so we also have a finiteness result for bi-arc-analytic homeomorphisms. On the other hand, the multiplicity of $f$
at the origin is given by
\[
\operatorname{mult}_0(f)
=
\min_\gamma\operatorname{ord}_0(f\circ\gamma),
\]
where $\gamma:(\mathbb R,0)\to(\mathbb C^n,0)$
is a real-analytic arc.
We use the convention
$\operatorname{mult}_0(0)=\infty$.
Thus, we obtain the following corollary:

\begin{corollary}
    Let $S_j$ be as in Theorem 1.1. Then, for all $t, s \in S_j$, it holds that $\text{mult}_0(f_t) = \text{mult}_0(f_s)$.
\end{corollary}

\begin{definition}
Let $f:(\mathbb C^n,0)\to(\mathbb C,0)$
be a holomorphic function germ. We define
\[
A(f)=
\left\{
\operatorname{ord}_0(f\circ\gamma):
\begin{array}{l}
\gamma:(\mathbb R,0)\to(\mathbb C^n,0)\\
\text{is a real-analytic arc}
\end{array}
\right\},
\]
where $\operatorname{ord}_0(0)=\infty$.
\end{definition}
The set $A(f)$ is known as the Fukui invariant
of $f$ (see \cite{IzumiKoikeKuo2002}).
\begin{corollary}
For any positive integers $n$ and $d$, the set
\[
\{A(f):f\in P_0(n,1,d;\mathbb C)\}
\]
is finite.
\end{corollary}

\begin{proof}
If $f=g\circ h$, where $h$ is a bi-arc-analytic
homeomorphism, then $A(f)=A(g)$, since composition
with $h$ induces a bijection between analytic arcs
through the origin and preserves the order of
$f$ along them. The result follows from the
finiteness of the bi-arc-analytic $\mathcal R$-equivalence
classes.
\end{proof}

Up to this point, we have considered polynomials
vanishing at the origin. We now turn to the case
of polynomials that do not necessarily vanish
at the origin, following the approach of
\cite[Section 4.3]{KoikePaunescu2022}. 
Let $P(n,1,d;\mathbb C)=\{f:\mathbb{C}^n \to \mathbb{C}; f \in \mathbb{C}[x_1, \cdots , x_n], deg(f) \leq d \}$.  In this setting, the definition of local
semialgebraic $\mathcal{RL}$-equivalence remains
the same, except that the homeomorphism $k$ in
the target maps $f(0)$ to $g(0)$ instead of
necessarily fixing the origin.

\begin{corollary}
For any positive integers $n$ and $d$, the number
of local semialgebraic $\mathcal{RL}$-equivalence
classes in $P(n,1,d;\mathbb C)$ is finite.

\end{corollary}

\begin{proof}
For each $f\in P(n,1,d;\mathbb C)$, set
$\widetilde f=f-f(0)$. By Theorem~1.1, there
are only finitely many local semialgebraic
$\mathcal R$-equivalence classes in
$P_0(n,1,d;\mathbb C)$. Moreover, the map
\[
\tau_f:\mathbb C\longrightarrow\mathbb C,
\qquad \tau_f(w)=w-f(0),
\]
is a semialgebraic homeomorphism satisfying
$\tau_f\circ f=\widetilde f$.
Thus $f$ and $\widetilde f$ are locally
semialgebraically $\mathcal{RL}$-equivalent. Therefore, every polynomial in $P(n,1,d;\mathbb{C})$ is $\mathcal{RL}$-equivalent to one of the finitely many representatives $g_1, \dots, g_r$ of $P_0(n,1,d;\mathbb{C})$, which proves the result.

\end{proof}
\begin{remark}
The analogous finiteness result does not hold
for local semialgebraic $\mathcal R$-equivalence
in $P(n,1,d;\mathbb C)$.
Indeed, the constant polynomials $f_c(z)=c$,
with $c\in\mathbb C$, are pairwise
$\mathcal R$-inequivalent.
\end{remark}

We now present one final application of Theorem~1.1.
In \cite{FernandesSampaioSilva2026}, Fernandes,
Sampaio, and Silva introduce the following notions:

\begin{definition}
Let $\alpha\in(0,1)$. Two polynomial functions
$f,g\in P_0(n,1,d;\mathbb C)$ are said to be
locally bi-$\alpha$-H\"older conjugate if there
exist neighborhoods $U,V$ of the origin and a
homeomorphism $h:U\to V$, with $h(0)=0$, such that
\[
f=g\circ h
\]
and both $h$ and $h^{-1}$ are $\alpha$-H\"older
continuous.

We say that $f$ and $g$ are locally H\"older
equivalent if they are locally bi-$\alpha$-H\"older
conjugate for every $\alpha\in(0,1)$.
In this case, we write $f\sim_H g$.
\end{definition}

They show that there exists $I \subset \mathbb C$ uncountable and a  family
$\{f_t\}_{t\in I}$ of complex polynomials in two
variables with uniformly bounded degree such that
$f_t$ and $f_s$ are not H\"older equivalent whenever
$t\neq s$ (see Remark~6.12).
In contrast, we introduce the following weaker
notion of equivalence:
\begin{definition}
Two polynomial functions
$f,g\in P_0(n,1,d;\mathbb C)$ are said to be
\emph{locally weakly H\"older equivalent}
if they are locally bi-$\alpha$-H\"older
conjugate for some $\alpha\in(0,1)$.
In this case, we write $f\sim_{wH}g$.
\end{definition}

This introduces an equivalence relation on $P_0(n,1,d;\mathbb C)$, since the composition of bi-$\alpha$-Hölder homeomorphisms with a bi-$\beta$-Hölder one results in a bi-$\alpha \beta$-Hölder homeomorphism. As germs of semialgebraic homeomorphisms are bi-$\alpha$-Hölder for some $\alpha$ by the {\L}ojasiewicz inequality, we have the following corollary:

\begin{corollary}
    The number of local weak H\"older equivalence classes in
$P_0(n,1,d;\mathbb C)$ is finite.
\end{corollary}
\begin{remark}
Consider the family of weighted homogeneous
polynomials $f_t(x,y)=x^3-3t^2xy^{2d}+y^{3d}, t\in\mathbb C$, where $d>1$ is an integer.
Fernandes, Sampaio and Silva
\cite{FernandesSampaioSilva2026}
showed that this family admits moduli under
local H\"older equivalence. In contrast, by Corollary~6.3, there exists a finite
subset $E\subset\mathbb C$ such that $f_t$ and $f_s$
are locally weakly H\"older equivalent,
for every $t,s\in\mathbb C\setminus E$.

\end{remark}

\section*{Declaration of AI use}

During the development of this work, the author
used ChatGPT (OpenAI) as an AI-assisted tool for
mathematical discussions, exploration of ideas,
and critical examination and refinement of arguments.
The author reviewed and independently verified
the mathematical statements and proofs, and assumes
full responsibility for the originality, accuracy,
and integrity of the manuscript.

% 

%==================================================
\section*{Acknowledgments}
% ==================================================

The author is grateful to Alexandre Fernandes
for his valuable mathematical discussions
and suggestions during the development
of this work.

% ==================================================
% BIBLIOGRAPHY
% ==================================================


\begin{thebibliography}{99}

\bibitem{BenedettiShiota1991}
R. Benedetti and M. Shiota,
Finiteness of semialgebraic types of polynomial functions,
\textit{Math. Z.}
\textbf{208} (1991), 589--596.

\bibitem{BochnakCosteRoy1998}
J. Bochnak, M. Coste and M.-F. Roy,
\textit{Real Algebraic Geometry},
Ergebnisse der Mathematik und ihrer Grenzgebiete,
vol. 36, Springer, Berlin, 1998.

\bibitem{FernandesSampaioSilva2026}
A. Fernandes, J. E. Sampaio and J. P. da Silva,
H\"older invariance of the Henry--Parusi\'nski invariant,
\textit{Mosc. Math. J.}
\textbf{26} (2026), no. 2, 145--161.

\bibitem{Fukuda1976Deformations}
T. Fukuda,
Deformations of real analytic functions and the
natural stratification of the space of real analytic functions,
\textit{Nagoya Math. J.}
\textbf{63} (1976), 139--152.

\bibitem{Fukuda1976Types}
T. Fukuda,
Types topologiques des polyn\^omes,
\textit{Publ. Math. Inst. Hautes \'Etudes Sci.}
\textbf{46} (1976), 87--106.

\bibitem{Fukuda1978}
T. Fukuda,
Lectures at Graduate School of Chiba University
(1978).

\bibitem{IzumiKoikeKuo2002}
S. Izumi, S. Koike and T.-C. Kuo,
Computations and stability of the Fukui invariant,
\textit{Compositio Math.}
\textbf{130} (2002), 49--73.

\bibitem{JansonDiscriminant}
S. Janson,
Resultant and discriminant of polynomials,
unpublished notes, revised 2010.
\url{https://people.kth.se/~jockum/Janson_Diskriminanter.pdf}.

\bibitem{KoikePaunescu2022}
S. Koike and L. P\u{a}unescu,
On finiteness theorems of polynomial functions,
\textit{European J. Math.}
\textbf{8} (2022), 214--234.

\bibitem{Milne2024}
J. S. Milne,
\textit{Algebraic Geometry},
version 6.10, course notes, 2024.
\url{https://www.jmilne.org/math/CourseNotes/AG.pdf}.

\bibitem{ParusinskiPaunescu2017}
A. Parusi\'nski and L. P\u{a}unescu,
Arc-wise analytic stratification, Whitney fibering
conjecture and Zariski equisingularity,
\textit{Adv. Math.}
\textbf{309} (2017), 254--305.

\bibitem{Shiota1987Nash}
M. Shiota,
\textit{Nash Manifolds},
Lecture Notes in Mathematics, vol. 1269,
Springer, Berlin, 1987.

\end{thebibliography}
\end{document}